\documentclass{article}

\usepackage{a4wide}

\usepackage{amsmath, amssymb, amsthm}
\usepackage{mathtools}
\usepackage{txfonts}
\usepackage{tikz}
\usepackage{pgfplots} 

\newtheorem{theorem}{Theorem}

\newtheorem{corollary}[theorem]{Corollary}
\newtheorem{definition}[theorem]{Definition}
\newtheorem{lemma}[theorem]{Lemma}
\newtheorem{proposition}[theorem]{Proposition}
\newtheorem{remark}[theorem]{Remark}

\newcommand{\eps}{\varepsilon}
\newcommand{\norm}[1]{\left\|#1\right\|}
\newcommand{\ord}[1]{\mathcal{O}\left(#1\right)}
\newcommand{\scal}[2]{\left\langle #1 , #2 \right\rangle}
\newcommand{\enorm}[1]{\left|\left|\left|#1\right|\right|\right|}
\newcommand{\lay}{BL}
\newcommand{\Olay}{{\Omega_{\lay}}}
\newcommand{\barOlay}{\bar{\Omega}_{\lay}}
\newcommand{\reg}{REG}
\newcommand{\Oreg}{{\Omega_{\reg}}}

\newcommand{\abs}[1]{\left|#1\right|}
\newcommand{\NN}{\mathbb{N}}
\newcommand{\RR}{\mathbb{R}}
\newcommand{\E}{\mathrm{e}}
\newcommand{\m}[1]{\mathcal{#1}}
\newcommand{\mB}{\mathcal{B}}
\newcommand{\mP}{\mathcal{P}}
\newcommand{\mS}{\mathcal{S}}
\newcommand{\mF}{\mathcal{F}}
\newcommand{\pt}{\partial}
\newcommand{\D}{\;\mathrm{d}}

\title{Balanced norm estimates for $rp$-Finite Element Methods applied to
  singularly perturbed fourth order boundary value problems}

\author{%
   Torsten Lin{\ss}\thanks{FernUniversit\"at in Hagen,
      Fakult\"at f\"ur Mathematik und Informatik,
      Universit\"atsstra{\ss}e 1, 58084, Hagen,
      email: \texttt{torsten.linss@fernuni-hagen.de}}
 \and Christos Xenophontos\thanks{University of Cyprus,
      Department of Mathematics \& Statistics, PO BOX 20537,
      1678, Nicosia,
      email: \texttt{xenophontos.christos@ucy.ac.cy}}
}

\begin{document}

\maketitle

\begin{abstract}
  We establish robust exponential convergence for $r p$--Finite Element Methods
  (FEMs) applied to fourth order singularly perturbed boundary value problems,
  in a \emph{balanced norm} which is stronger than the usual energy norm
  associated with the problem.   As a corollary, we get robust exponential convergence in the maximum norm.
  $r p$ FEMs are simply $p$ FEMs with possible repositioning of the (fixed number of) nodes.
  This is done for a $C^1$ Galerkin FEM in 1-D, and a $C^0$ mixed FEM in 2-D
  over domains with smooth boundary. In both cases we utilize the
  \emph{Spectral Boundary Layer} mesh from \cite{MXO}.


\end{abstract}

%%%%%%%%%%%%%%%%%%%%%%%%%%%%%%%%%%%%%%%%%%%%%%%%%%%%%%%%%%%%%%%%%

\section{Introduction}
\label{intro}
The quest for \emph{balanced norm} estimates for singularly perturbed problems (SPPs) 
has been going on for a while now (see \cite{FR}--\cite{LS} and the references therein). 
When using the Finite Element Method (FEM),
the error estimates are usually given in the \emph{energy norm} defined by the bilinear form in the
variational formulation of the problem. In certain cases, this norm is too weak and cannot see 
the layers present in the solution (see, e.g., \cite{LS}). While it is straight forward to define a balanced norm,
the bilinear form lacks coercivity with respect to it, hence the usual methods of proof cannot be applied. Various
approaches have been proposed \cite{FR}--\cite{LS}, dealing mainly with reaction-diffusion type problems.
The majority of these works pertains to fixed order FEMs; notable exceptions are \cite{CVX, MX}. 

In this article we consider fourth order SPPs and the approximation 
of their solution by the $r p$ version of the FEM on the \emph{Spectral Boundary Layer} mesh from \cite{MXO}.
In this version the position (but not the number) of the nodes could possibly change, in addition to 
the polynomial degree $p$ of the approximating polynomials, being increased.
In the literature, one finds results on robust exponential convergence in the energy norm \cite{PC, CFLX, PZMX}.
We prove robust exponential convergence in a balanced norm and as a corollary we get robust exponential
convergence in the maximum norm. In one-dimension this is achieved for a $C^1$ Galerkin FEM, used in \cite{PC, PZMX},
and in two-dimensions for the $C^0$ mixed formulation from \cite{CFLX}. The novelty in our approach lies with 
the choice of a certain projection
of the solution onto the FE space. Specifically, we use an interpolant in the layer region and an appropriate
projection in the rest of the domain. While this general approach has been used before (e.g., \cite{FR3}),
our contribution differs twofold: first we consider the $rp$ version (hence exponential rates of convergence 
are obtained), and second, the projection is used \emph{only} for the smooth parts of the solution, as opposed
to the entire solution as was done in \cite{FR3}. Finally, our contribution proves the conjectures made in \cite{PC} and \cite{CFLX} based on numerical evidence.

The rest of the paper is organized as follows: in Section \ref{sec:1D} we consider a one-dimensional problem
and its discretization by a $C^1$ Galerkin FEM. We present the available estimates 
from the literature and in Section \ref{sec:balanced} we improve them by establishing 
robust exponential convergence in a balanced, as well as the maximum norm.
Section \ref{sec:2D} contains the two-dimensional case, over domains with smooth boundary,
and the $C^0$ mixed method from \cite{CFLX}. We, again, describe the estimates from the literature 
and then in Section \ref{estimates} we improve
them, just like in one-dimension.

\paragraph{Notation}
With \mbox{$D \subset \RR^d$}, \mbox{$d=1,2$}, a domain with boundary $\pt D$
and measure $\abs{D}$, we will denote by $C^{k}(D)$ the space of continuous
functions on $D$ with continuous derivatives up to order~$k$.
We will use the usual Sobolev spaces $W^{k,q}(D)$, \mbox{$q\in[0,\infty]$},
\mbox{$k\in\NN_0$} of functions on $D$
generalized derivatives of order $0,1,\dots,k$ in $L^{q}(D)$, equipped with
the norm and seminorm $\norm{\cdot}_{k,q,D}$ and $\abs{\cdot}_{k,m,D}$,
respectively.
When $q=2$, we will write $H^k(D)$ instead of $W^{k,2}(D)$, and for the norm
and seminorm, we will write $\norm{\cdot}_{k,D}$ and $\abs{\cdot}_{k,D}$,
respectively.
The usual $L^{2}(D)$ inner product will be denoted by $\scal{\cdot}{\cdot}_{D}$,
with the subscript omitted when \mbox{$D=\Omega$} and there is no confusion.
We will also use the space 
\begin{gather*}
  H_0^2\left(D\right) = \left\{ u\in H^2(D) \colon
      \left.u\right\vert _{\pt D} = \left.\frac{\pt u}{\pt n}\right\vert _{\pt D}=0\right\},
\end{gather*}
where $\frac{\pt}{\pt n} $ denotes the normal derivative. 
The norm of the space $L^\infty(D) = W^{0,\infty}(D)$
of essentially bounded functions on $D$ is denoted by
$\norm{\cdot}_{L^\infty(D)} = 
  \norm{\cdot}_{0,\infty,D} = \abs{\cdot}_{0,\infty,D}$.

Finally, the notation ``$a\lesssim b$'' means ``$a\leq Cb$''
with $C$ being a generic positive constant, independent of any
discretization or singular perturbation parameters.

\section{One-dimensional problems}
\label{sec:1D}
Since the boundary layer effect is one-dimensional (in the direction normal to the boundary), studying
one-dimensional problems is a necessary first step towards two (and higher) dimensions. In fact, most of the
ideas are present in the one-dimensional case, and their extention to two-dimensions is (in most cases) 
straight forward.

We consider the following problem (cf. \cite{PC, PZMX}): find $u$ such that
\begin{subequations}\label{de-bc}
  \begin{alignat}{2}
    \eps^{2}u^{(4)}-\left( bu'\right)'+cu & = f  &\quad& \text{in}
       \ \ \Omega \coloneqq \left(0,1\right),  \\ %\label{de} \\
      u(0)=u(1)=u'(0)=u'(1) & = 0,  \label{bc}
  \end{alignat}
\end{subequations}
where \mbox{$0<\eps\leq 1$} is a given parameter
that can approach zero and the functions $b$, $c$ and $f$ are given and
sufficiently smooth.
In particular, we assume that they are (real) \emph{analytic} functions
satisfying, for some positive constant $\gamma_d$, independent of $\eps$, 
\begin{gather}
  \label{analytic}
  \abs{f}_{n,\infty,\Omega}
   + \abs{c}_{n,\infty,\Omega}
   + \abs{b}_{n,\infty,\Omega}\lesssim n!\gamma_d^n
    \quad \forall \; n\in \NN_0.
\end{gather}
In addition, we assume that there exist positive constants
$b_{\min}$ and $c_{\min}$,
independent of $\eps$, such that
\begin{gather}
  \label{data}
  b\ge b_{\min}^2, \quad c\geq c_{\min}^2 \quad \text{on} \ \
     \bar{\Omega}, 
\end{gather}

It is well known (see, e.g. \cite{OMalley}) that the solution $u$ to
\eqref{de-bc} can be decomposed into a smooth part, two boundary layer parts
and a remainder.
We have the following result from~\cite{PZMX}.

\begin{proposition}
  \label{prop:reg}
  Let $u$ be the solution to \eqref{de-bc}, and assume \eqref{analytic} holds.
  Then, there exist a positive constants $K$, $\gamma$ and $\delta$ such that
  for all \mbox{$n\in \NN_0$}, 
  \begin{gather} \label{A1}
    \abs{u}_{n,\infty,\Omega} \lesssim K^{n}\max \left\{ n^{n}, \varepsilon^{1-n}\right\}\,.
  \end{gather}
  Moreover, $u$ may be decomposed as
  \begin{subequations}
    \label{AA}
    \begin{gather}
      % \label{A2}
        u=u_S+\tilde{u} + \bar{u} + u_R,
    \end{gather}%
    with
    \begin{alignat}{2}
      % \label{A3}%
        \abs{u_S}_{n,\infty,\Omega} & \lesssim K^nn^n, \\
      \label{A4}%
        \abs{\tilde{u}^{(n)}(x)} + \abs{\bar{u}^{(n)}(1-x)}
                   & \lesssim K^{n} \eps^{1-n} 
                     \E^{-\gamma x/\eps}\ \ \forall \ x\in\bar\Omega \\
      %\intertext{and}
        \label{A6}
        \norm{u_R}_{2,\Omega} & \lesssim \E^{-\delta /\eps}.
    \end{alignat}
  \end{subequations}
  In this decomposition, $u_S$ denotes the smooth part,
  $\tilde{u}$ the boundary layer at the left endpoint,
  $\bar{u}$ denotes the  boundary layer at the right endpoint, and $u_R$ is the remainder.
\end{proposition}

We mention that \eqref{A1} corresponds to classical differentiability, while
\eqref{AA} correspond to regularity through asymptotic expansions, see \cite{PC,melenk}.

%%%%%%%%%%%%%%%%%%%%%%%%%%%%%%%%%%%%%%%%%
\subsection{Discretization by a $C^{1}$ $rp$-FEM}% \label{mesh}}

The variational formulation of \eqref{de-bc} reads:
Find $u\in H_{0}^{2}(\Omega)$ such that 
\begin{gather}
  \label{BuvFv}
    \mB_\eps \left(u,v\right) = \mF \left(v\right)
      \;\;\forall \;v\in H_{0}^{2}(\Omega),
  \intertext{where}
  \label{Buv}
    \mB_\eps \left(w,v\right) \coloneqq
      \eps^2 \scal{w''}{ v''}_{\Omega}
       + \scal{bw'}{v'}_{\Omega} + \scal{cw}{v}_{\Omega},
       \quad
    \mF(v) \coloneqq \scal{f}{v}_{\Omega}.
\end{gather}
Because of~\eqref{data}, the bilinear form $\mB_{\eps}\left(\cdot,\cdot\right)$
induces a norm
\begin{gather}
  \label{energy}
  \norm{w}_{E} \coloneqq \mB_{\eps}\left(w,w\right) ,
\end{gather}%
the so called \emph{energy norm}.
Clearly, $\mB_{\eps}\left(\cdot,\cdot\right)$ is coercive (with constant $1$)
with respect to this norm.

We shall seek an approximation to $u$ in a finite dimensional subspace
of~$H_0^2(\Omega)$.
To this end, let $\mP_p$, \mbox{$p\in\NN_0$}, denote the space of polynomials of
maximum degree $p$, and let \mbox{$\Delta \coloneqq \left\{x_i\right\} _{i=0}^{N}$}
be an arbitrary subdivision of $\Omega$ with mesh intervals \mbox{$\Omega_j\coloneqq(x_{i-1},x_i)$}
and mesh sizes \mbox{$h_i\coloneqq x_i-x_{i-1}>0$}, \mbox{$i=1,\dots,N$}.
We define the finite dimensional spaces
\begin{gather*}
  \mS^p_\Delta \coloneqq
     \Bigl\{w\in C^1(\bar\Omega) \colon w\vert_{\Omega_i}\in\mP_p, \
        i=1,\dots,N\Bigr\}\,, \quad
  \mS_{0,\Delta}^p \coloneqq \mS^p_\Delta \cap H_0^2(\Omega)\,.
\end{gather*}

Our discretization of \eqref{de-bc} reads: 
Find \mbox{$u_{p,\Delta}\in \mS_{0,\Delta}^p$} such that 
\begin{gather*}
  \mB_\eps\left(u_{p,\Delta},v\right) = \mF(v) \;\;\forall \;v\in \mS_{0,\Delta}^p.
\end{gather*}
This discretization possesses the Galerkin orthogonality property
\begin{gather*}
  \mB_{\eps} \left( u-u_{p,\Delta},v\right) = 0\;\;\forall \;v\in  \mS_{0,\Delta}^p
\end{gather*}
as well as best approximation property
\begin{gather}
  \label{bestapprox}
  \norm{u-u_{p,\Delta}}_{E} = \inf_{v\in  \mS_{0,\Delta}^p}
      \norm{u-v}_{E}\,.
\end{gather}%

Given \mbox{$w\in H^2(\Omega)$}, an interpolant
\mbox{$I^p_\Delta w \in \mS^p_\Delta$}, \mbox{$p\ge 3$},
is uniquely defined by
\begin{align*}
    \left(I^p_\Delta w - w\right)^{(k)} (x_i) & = 0, \quad i=0,\dots,N, \ \ k=0,1 \\
  \intertext{and}
    \int_{\Omega_j}\left(I^p_\Delta w - w\right)'' (x) q(x) \D x & = 0, \ \forall q \in \mP_{p-2}, \ \ i=0,\dots,N.
\end{align*}

\begin{lemma}
  Let $w\in H^\ell(\Omega_j)$, $\ell\ge 2$.
  Then for all $s$, $0\le s \le \min(p-1,\ell-2)$,
  the interpolation error satisfies the local bounds
  \begin{gather*}
    (p-1)^{2-k}\abs{w-I_\Delta^p w}_{k,\Omega_j}^2 \le
       \left(\frac{\abs{\Omega_j}}{2}\right)^{2(s+2-k)}
       \frac{\left(p-1-s\right)!}{\left(p+1+s\right)!}
       \abs{w}_{s+2,\Omega_j}^2\,, \ \ k=0,1,2.
  \end{gather*}
\end{lemma}
\begin{proof}
  See, \cite{schwab} and \cite{BdVBRS}.
\end{proof}

% An interpolant $I_p$ may be defined on the above mesh, as the truncated Legendre series of $u''$ (see, e.g. \cite{PZMX}).
% The approximation properties of $I_p$ are given below.

We next give the definition of the \emph{Spectral Boundary Layer Mesh}
we will utilize throughout the article.
\begin{definition}% [Spectral Boundary Layer mesh]
  \label{SBL}
  For \mbox{$\kappa>0$}, \mbox{$p\in \NN$} and \mbox{$0 < \eps \le 1$},
  define the Spectral Boundary Layer mesh $\Delta _{BL}$ as
  \begin{gather*}
    \Delta _{BL} \coloneqq \bigl\{0,\tau,1-\tau,1\bigr\}\,,\quad
      \tau \coloneqq \min\left\{\kappa p\eps, 1/3\right\}\,,
  \end{gather*}
  with the layer and coarse-mesh regions

  \begin{gather*}
    \Omega_\ell \coloneqq (0, \tau) \cup (1-\tau,1) \, , \ \
    \Omega_c    \coloneqq (\tau,1-\tau).
  \end{gather*}
\end{definition}

We shall consider the sequence \mbox{$u_p\in\mS^p_{0,\Delta_{BL}}$},
\mbox{$p=3,4,\dots$}, of approximations defined by
\begin{gather}
  \label{discrete}
  \mB_\eps\left(u_p,v\right) = \mF(v) \;\;\forall \;v\in \mS_{0,\Delta_{BL}}^p.
\end{gather}
Given \mbox{$w\in H^2(\Omega)$}, we define interpolants
\mbox{$I_p w \coloneqq I^p_{\Delta_{BL}} w$}, \mbox{$p=3,4,\dots$}

The best approximation property~\eqref{bestapprox} implies
\begin{gather*}
  \norm{u-u_p}_{E} \le \norm{u-I_p u}_{E}\,, \ \ p=3,4,\dots
\end{gather*}

\begin{proposition}\label{prop:interp}
  Let $u$ be the solution to \eqref{de-bc} and assume that \eqref{analytic}
  holds.
  Then there exist positive constants $\kappa_0$, $\kappa _1$ and $\beta>0$
  independent of $\eps$ and $p$, such that the following is true:
  For every \mbox{$p\in\NN$} and every \mbox{$\kappa \in (0,\kappa _{0}]$}
  with \mbox{$\kappa p \geq \kappa_1$}, the interpolation error satisfies
  \begin{equation}\label{interp_est}
    \norm{u-I_{p}u} _{1,\infty,\Omega}
      + \eps^{1/2} \abs{u-I_p u}_{2,\Omega} \lesssim \E^{-\beta p}.
  \end{equation}
 Moreover,
  \begin{equation}\label{interp_uS}
    \norm{u_S-I_{p}u_S} _{2,\Omega} \lesssim \E^{-\beta p},
  \end{equation}
and
  \begin{equation}\label{interp_lay}
    \abs{u- I_p u}_{k,\Olay} \lesssim \eps^{3/2-k} \E^{-\beta p}.
  \end{equation}

%  Moreover, in the layer and coarse-mesh regions there holds\footnote{All needed?}
%%
%\begin{eqnarray}
%\eps^{-3/2} \Vert u- I_p u \Vert_{L^{\infty}(\Olay)} +\eps^{-1/2} \Vert \left( u- I_p u \right)' \Vert_{L^{\infty}(\Olay)} +\eps^{-1} \Vert u- I_p u \Vert_{0,\Olay}  \notag \\
%+ \vert u- I_p u \vert_{1,\Olay} +
%\xi^{1/2} \vert v- I_p v \vert_{2,\Olay} \lesssim C_u  e^{-\beta p} \label{res1}, \\
%%
%\varepsilon^{-1}\Vert u- I_p u\Vert_{L^{\infty}(\Oreg)} + \Vert \left(u- I_p u\right)' \Vert_{L^{\infty}(\Oreg)}
%+\varepsilon^{-3/2} \Vert u- I_p u \Vert_{0,\Oreg} \notag \\
%+\varepsilon^{-1/2}\vert u- I_p u \vert_{1,\Oreg}
%+ \varepsilon^{1/2} \vert u- I_p u \vert_{2,\Oreg} \lesssim C_u e^{-\kappa p}. \label{res2}
%\end{eqnarray}
\end{proposition}
\begin{proof}
  Inequality \eqref{interp_est} was established in \cite{PZMX} (see also \cite{PC}),
  and inequalities \eqref{interp_uS}, \eqref{interp_lay} are shown in the Appendix.
\end{proof}

Using \eqref{bestapprox} and~\eqref{interp_est} the following may be established, cf.~\cite{PC}.
\begin{proposition}
  \label{thm:main}
  Let $u$ solve \eqref{de-bc} and let $u_p$, $p=3,4,\dots$, be its finite
  element approximations obtained by~\eqref{discrete}.
  Then there exists a constant $\sigma >0$, independent of $\eps$ but depending on the data, such that
  \begin{gather*}
    \norm{u-u_p}_{E} \lesssim \E^{-\sigma p}.
  \end{gather*}
\end{proposition}

%%%%%%%%%%%%%%%%%%%%%%%%%%%%%%%%%%%%%%%%%%%%%%%%
\subsection{Balanced norm estimates}
\label{sec:balanced}
The energy norm used in Proposition \ref{thm:main} is not the appropriate norm for this problem, since with 
\begin{gather*}
  u_{BL} \coloneqq \tilde{u}+ \bar{u},
\end{gather*}
we calculate
\begin{gather*}
  \norm{u_{BL}}_{E} = \ord{\eps^{1/2}} \; \mbox{ while } \; 
  \norm{u_{S}}_{E}  = \ord{1}\,, \ \ \eps \to 0.
\end{gather*}
Thus, the energy norm does not see the layers present in the solution
of \eqref{de-bc}, as $\eps \to 0$. 
The correct weight for $\abs{u}_{2,\Omega}$ is $\eps^{1/2}$
(as opposed to $\eps$, cf. \eqref{energy}).
In which case, the resulting norm would be \emph{balanced}, like the one below:
\begin{equation}\label{balanced-1D}
  \norm{u}^2_{B} \coloneqq \eps \abs{u}_{2,\Omega}^2 + \norm{u}_{1,\Omega}^2.
\end{equation}
In this norm we have
\begin{gather*}
  \norm{u_{BL}}_{B} = \ord{1} = \norm{u_{S}}_{B}\,, \ \ \eps \to 0.
\end{gather*}
Unfortunately, the bilinear form $\mathcal{B}_{\eps}$ given by \eqref{Buv},
is not coercive with respect to this norm.
Nevertheless, numerical experiments reported in \cite{PC}, suggest the following:
\begin{gather}
  \label{C1}% \label{C3}
  \norm{u - u_p}_{1,\infty,\Omega} + \norm{u-u_p}_B \lesssim \E^{-\sigma p}
\end{gather}
Our goal in this section is to establish \eqref{C1}.

We begin by studying the auxiliary problem of finding
\mbox{$u_{S,p} \in \mP_p$}, \mbox{$p\ge3$}, such that
\begin{align} \label{pi_p}
  \mB_{0}\left(u_{S,p}-u_{S},v\right) = 0
    \ \ \forall \; v \in \mP_p\cap H_0^1(\Omega), \ \
               u_{S,p}-u_{S} = 0 \ \text{on} \ \pt\Omega
\end{align}
where
\begin{gather*}
  \mB_{0}\left(w,v\right) \coloneqq \scal{b w'}{v'} + \scal{c w}{v}.
\end{gather*}
Note that $\mB_{0}$ induces a norm
\mbox{$\enorm{w} \coloneqq \mB_{0}(w,w)^{1/2}$}
that is equivalent to the standard
$H^1(\Oreg)$ norm because of~\eqref{data}.
Consequently, $u_{S,p}$ is uniquely defined and we have the following.
\begin{lemma}\label{lem:pi_p}
  Let \mbox{$u_{S,p} \in \mP_p$}, \mbox{$p\ge3$}, be defined by \eqref{pi_p}.
  Then, 
  \begin{gather*}
    \abs{u_{S,p} - u_S}_{k} \lesssim \E^{-\beta p}, \quad k=0,1,2, \\
    \norm{u_{S,p} - u_S}_{C^1(\bar\Omega)} \lesssim \E^{-\beta p}
    \intertext{and}
    \abs{\left(u_{S,p} - u_S\right)(x)} \lesssim \min(x,1-x) \E^{-\beta p}\,,
      \ \ x\in \bar\Omega.
  \end{gather*}
\end{lemma}
\begin{proof}
  Note, that the bilinear form in~\eqref{pi_p} is coercive and bounded.
  So standard techniques for spectral methods apply to give the bounds of
  the first inequality for $k=0$ and $k=1$.

  For $k=2$ we proceed as follows.
  We have
  \begin{gather*}
    \enorm{u_{S,p} - u_S} \le \enorm{v - u_S}
       \ \ \forall \ v\in\mP_p \ \ \text{with} \ \ v= u_S \ \ \text{on} \ \pt \Omega,
  \end{gather*}
  which implies
  \begin{gather*}
    \enorm{u_{S,p} - I_p u_S} \le 2 \enorm{I_p u_S - u_S},
  \end{gather*}
  by the triangle inequality.
  Next
  \begin{align*}
    \abs{u_{S,p} - u_S}_{2}
      & \le \abs{u_{S,p} - I_p u_S}_{2}
          + \abs{I_p u_S - u}_{2}
        \le (p-1)^2 \abs{u_{S,p} - I_p u_S}_{1}
          + \abs{I_p u_S - u}_{2} \\
      & \le 2 b_0^{-1} (p-1)^2 \enorm{I_p u_S-u_S}_{1}
          + \abs{I_p u_S - u}_{2}\,,
  \end{align*}
  where we have used an inverse inequality, see~Lemma~\ref{lem:inverse-1D}.
  Now the desired result follows from~\eqref{interp_uS}.

  A Sobolev embedding gives the bound in the $C^1$ norm, while the final bound
  follows upon noting that $\left(u_{S,p}-u_S\right)(x)=0$ for $x\in\{0,1\}$
  and the maximum-norm bound for the derivative.
\end{proof}

Next, we construct a special FE-function \mbox{$\hat{u}_p\in \mS^p_{0,\Delta_{BL}}$}
that is particularly close to $u$.
To this end, define auxiliary functions \mbox{$\chi_0,\chi_1\in \mP_3$} by
\begin{gather*}
  \chi_0(0)=\chi_0'(0)=\chi_0(\tau)=0, \quad \chi_0'(\tau)=1, \\
 \intertext{and}
  \chi_1(0)=\chi_1'(0)=\chi_1'(\tau)=0, \quad \chi_1(\tau)=1.
\end{gather*}
A direct calculation establishes the following bounds on the $\chi_i$:
\begin{gather}\label{bound-chi}
  \abs{\chi_i}_{k,(0,\tau)}
     \lesssim \tau^{3/2-k-i}, \; k \in \{0, 1, 2\}, \; i \in \{0,1\}.
\end{gather}
Now, set
\begin{gather}\label{uhat}
  \hat{u}_p(x) \coloneqq
     \begin{cases}
        \left(I_p u\right)(x) + \chi_0(x)\left(u_{S,p} - u\right)'(\tau) % \\
      %\qquad\qquad\quad
           + \chi_1(x)\left(u_{S,p} - u\right)(\tau), & x \in [0, \tau], \\
        u_{S,p}(x), & x \in [\tau, 1-\tau], \\
      (I_p u)(x)-\chi_0(1-x)\left(u_{S,p} - u\right)'(1-\tau), % \\
      % \qquad\qquad\quad
           + \chi_1(1-x)\left(u_{S,p}-u\right)(1-\tau), & x \in[1-\tau,1].
     \end{cases}
\end{gather}

We re-iterate that in the layer region, $\hat{u}_p$ is given as the
interpolant $I_pu$ of $u$ plus corrections that ensure
$\hat{u}_p\in C^1(\bar\Omega)$.
In the regular region, we \emph{only} apply the projection on $u_S$.
We have the following.
\begin{lemma}\label{lem:u-ut}
  Assume that~\eqref{AA} holds and that $\hat{u}_p$ is defined
  by \eqref{uhat}.
  Then the following bounds hold true:
  \begin{subequations}
  \begin{alignat}{2} \label{u-ut:bl}
    \abs{u - \hat{u}_p}_{k, \Olay}
       & \lesssim \eps^{3/2-k} \E^{-\beta p}, & \ & k \in \{0,1,2\}, \\
    \label{u-ut:reg}
    \abs{u - \hat{u}_p}_{k, \Oreg}
       & \lesssim \max\left(1,\eps^{3/2-k}\right)
            \E^{-\beta p}, & \ &  k\in \{0,1,2\}.
  \end{alignat}
  \end{subequations}
\end{lemma}

\begin{proof}
  We shall study the two regions $\Olay$ and $\Oreg$ separately.
  \paragraph{$\Olay$:}
  We present the argument for $(0,\tau)$,
  because identical bounds hold for the interval $(1-\tau,1)$.
  We have
  \begin{gather*}
    \abs{u - \hat{u}_p}_{k,(0,\tau)}
        \le \abs{u - I_p u}_{k,(0,\tau)}
        + \abs{\chi_0}_{k,(0,\tau)}
            \cdot \abs{\left(u_{S,p} - u\right)'(\tau)}
        + \abs{\chi_1}_{k,(0,\tau)}
            \cdot \abs{\left(u_{S,p} - u\right)(\tau)},
        \ \ k\in\{0,1,2\}.
  \end{gather*}
  Proposition~\ref{prop:interp} %,~\eqref{interp_lay}
  takes care of the first term above.
  For the other two terms, a triangle inequality,
  Lemma \ref{lem:pi_p}, \eqref{A4} and \eqref{A6} yield
  \begin{gather*}
     \abs{\left(u_{S,p} - u\right)^{(\ell)}(\tau)}
       = \abs{\left(u_{S,p} - u_S - u_{BL} - u_R\right)^{(\ell)}(\tau)}
       \lesssim \tau^{1-\ell} \E^{-\beta p}\,, \ \ \ell\in{0,1}.
  \end{gather*}
  Recalling~\eqref{bound-chi}, we get
  \begin{gather*}
     \abs{\chi_{1-\ell}^{}}_{k,(0,\tau)}
     \abs{\left(u_{S,p} - u\right)^{(\ell)}(\tau)}
       \lesssim \tau^{3/2-k} \E^{-\beta p}\,, \ \ \ell\in{0,1},
  \end{gather*}
  and~\eqref{u-ut:bl} follows.
  We mention here that we will routinely ``hide'' polynomial powers into exponentials, e.g.
  $p^s \E^{-\beta p} \lesssim \E^{-\beta p} \; \; \forall \;\; s \in \mathbb{R}$.

  \paragraph{$\Oreg$:}
  We have from the definition of $\hat{u}_p$,
  \begin{gather*}
    \abs{u -\hat{u}_p}_{k,\Oreg} \le
      \abs{u_{S}-u_{S,p}}_{k,\Oreg}
        + \abs{u_{BL}}_{k,\Oreg}
        + \abs{u_R}_{k,\Oreg}, \ k=0,1,2.
  \end{gather*}
  The first term has been bounded in Lemma~\ref{lem:pi_p}.
  For $u_{BL}$ we have by~\eqref{A4}, and because $x\in[\tau,1-\tau]$,
  \begin{gather*}
    \abs{u_{BL}}_{k,\Oreg} \lesssim \eps^{3/2-k} \E^{-\gamma \kappa p}, \ k=0,1,2,
  \end{gather*}
  while for $u_R$, ineq.~\eqref{A6} gives
  \begin{gather*}
    \abs{u_R}_{k,\Oreg} \lesssim \E^{-3\delta \kappa p}, \ k=0,1,2.
  \end{gather*}
  Combining the last two inequalities with Lemma~\ref{lem:pi_p},
  we obtain~\eqref{u-ut:reg}.
\end{proof}

We now have the necessary tools for proving the following, main result of this section.

\begin{theorem}\label{thm:C1balanced}
  Let $u$ be the solution of \eqref{BuvFv} and let $u_p$ be the approximation
  obtained by \eqref{discrete}.
  Then, there exists a positive constant $\beta$,
  independent of $\eps$, such that
  \begin{gather*}
    \norm{u - u_p}_{B} \lesssim \E^{-\beta p}.
  \end{gather*}
\end{theorem}
\begin{proof}
  By~\eqref{balanced-1D}
  \begin{gather*}
    \norm{u-u_p}^2_{B} = \eps  \abs{u - u_p}^2_{2,\Omega}
      + \norm{u-u_p}^2_{1,\Omega}.
  \end{gather*}
  The only ``troublesome'' term above is $\eps \abs{u-u_p}_{2,\Omega}$.
  (The other one is handled by Proposition \ref{thm:main}).
  To deal with it, we use the triangle inequality (with $\hat{u}_p$ defined
  by \eqref{uhat}) to get
  \begin{gather*}
    \eps^{1/2} \abs{u - u_p}_{2,\Omega}
      \le \eps^{1/2} \abs{u-\hat{u}_p}_{2,\Omega}
               + \eps^{1/2} \abs{\hat{u}_p-u_p}_{2,\Omega}.
  \end{gather*}
  Only the term $\eps^{1/2} \abs{\hat{u}_p-u_p}_{2,\Omega}$ needs to be
  considered, since the rest can be bounded by Lemma \ref{lem:u-ut}.

  Let  $\eta \coloneqq \hat{u}_p-u_p \in \mS^p_{0,\Delta_{BL}}$.
  Galerkin orthogonality gives
  \begin{gather*}
    \norm{\eta}_E^2 =  \eps^2 \scal{\left(\hat{u}_p-u\right)''}{\eta''} +
        \scal{b(\hat{u}_p-u)'}{\eta'} + \scal{c(\hat{u}_p-u)}{\eta} .
  \end{gather*}
  On the coarse mesh region $\Omega_{REG}$, we have
  \begin{align*}
    & \scal{b\left(\hat{u}_p-u\right)'}{\eta'}_{\Oreg}
      + \scal{c\left(\hat{u}_p-u\right)}{\eta}_{\Oreg} \\
    & \qquad =
      \scal{b\left(u_{S,p}-u_S\right)'}{\eta'}_{\Oreg}
      + \scal{c\left(u_{S,p}-u_S\right)}{\eta}_{\Oreg} \\
    & \qquad\qquad\qquad
      - \scal{b\left(u_{BL}+u_R\right)'}{\eta'}_{\Oreg}
      - \scal{c\left(u_{BL}+u_R\right)}{\eta}_{\Oreg} \\
    & \qquad =
      - \scal{b\left(u_{S,p}-u_S\right)'}{\eta'}_{\Olay}
      - \scal{c\left(u_{S,p}-u_S\right)}{\eta}_{\Olay} \\
    & \qquad\qquad\qquad
      - \scal{b\left(u_{BL}+u_R\right)'}{\eta'}_{\Oreg}
      - \scal{c\left(u_{BL}+u_R\right)}{\eta}_{\Oreg},
  \end{align*}
  where the definition of $\hat{u}_p$ was used.
  Now, the $C^1$-norm bounds of Lemma \ref{lem:pi_p},
  \eqref{A4} and \eqref{A6} yield
  \begin{align*}
    & \abs{\scal{b\left(\hat{u}_p-u\right)'}{\eta'}_{\Oreg}
      + \scal{c\left(\hat{u}_p-u\right)}{\eta}_{\Oreg}} \\
    & \qquad
      \lesssim \tau^{1/2} \norm{u_{S,p}-u_S}_{C^1(\barOlay)} \norm{\eta}_{1,\Olay}
       + \norm{u_{BL}+u_R}_{1,\Oreg} \norm{\eta}_{1,\Oreg}
      \lesssim  \eps^{1/2} \E^{-\beta p} \norm{\eta}_{1,\Omega}.
  \end{align*}
  Similarly, using the $H^k$-norm bounds of Lemma \ref{lem:pi_p},
  \eqref{A4} and \eqref{A6}, we get
  \begin{align*}
    & \eps^2\abs{\scal{(\hat{u}_p-u)''}{\eta''}_{\Oreg}} \\
    & \qquad \le
      \eps^2 \abs{\scal{(u_{S,p}-u_S)''}{\eta''}_{\Oreg}}
        + \eps^2\abs{\scal{\left(u_{BL}+u_R\right)''}{\eta''}_{\Oreg}}
      \lesssim \eps^{3/2} \E^{-\beta p} \abs{\eta}_{2,\Oreg}.
  \end{align*}

  In the layer region $\Olay$, we use the Cauchy-Schwarz inequality
  and~\eqref{u-ut:bl} to obtain
  \begin{gather*}
    \abs{\scal{b\left(\hat{u}_p-u\right)'}{\eta'}_{\Olay}
      + \scal{c\left(\hat{u}_p-u\right)}{\eta}_{\Olay}}
      \lesssim \eps^{1/2} \E^{\beta p} \norm{\eta}_{1,\Olay}
    \intertext{and}
    \eps^2\abs{\scal{(\hat{u}_p-u)''}{\eta''}_{\Olay}}
      \lesssim \eps^{1/2} \E^{\beta p} \eps \abs{\eta}_{2,\Olay}.
  \end{gather*}

  Combining the above we have
  \begin{gather*}
    \norm{\eta}_E^2 \lesssim  \eps^{1/2} \E^{-\beta p} \norm{\eta}_E.
  \end{gather*}
  Dividing by $\norm{\eta}_E$ yields
  \begin{gather*}
    \eps^{1/2} \abs{u_p-\hat{u}_p}_{2,\Omega}
      \lesssim \norm{u_p-\hat{u}_p}_E \lesssim \E^{-\beta p},
  \end{gather*}
  and from the triangle inequality we have
  \begin{gather*}
     \eps^{1/2} \abs{u_p-u}_{2,\Omega} \lesssim \E^{-\beta p}.
  \end{gather*}
  This completes the proof.
\end{proof}

Once we have the balanced norm estimate, we may obtain a maximum norm estimate as follows.
\begin{corollary}
  \label{max_norm_C1}
  Let $u$ be the solution of \eqref{BuvFv} and let $u_p$ be the solution
  of~\eqref{discrete}, based on the Spectral Boundary Layer mesh of
  Definition~\ref{SBL}.
  Then, there exists a positive constant $\beta$, independent of $\eps$, such that
  \begin{gather*}
    \norm{u - u_p}_{C^1(\bar\Omega)} \lesssim \E^{-\beta p}.
  \end{gather*}
\end{corollary}
\begin{proof}
  Let $x \in [0,\tau]$ be arbitrary.
  Then, using the boundary conditions \eqref{bc}, we have
  \begin{gather*}
    (u-u_p)^{(k)}(x)
      = \int_0^{x} \left(u-u_p\right)^{(k+1)}(t) \D t , \  k \in \{0, 1\},
  \end{gather*}
  and by the Cauchy-Schwarz inequality
  \begin{gather*}
    \abs{(u-u_p)^{(k)}(x)}
      \lesssim \tau^{1/2} \abs{u-u_p}_{k+1,\Olay}
      \lesssim (p\eps)^{1/2} \eps^{-1/2} \E^{-\beta p}
      \lesssim \E^{-\beta p},
  \end{gather*}
  where Theorem \ref{thm:C1balanced} was used.

  The same technique works for the other layer region $[1-\tau, 1]$,
  so let us consider the coarse region $\Omega_{\mathrm{REG}}$.
  We use the triangle inequality, with $I_p$ the interpolant of
  Proposition \ref{prop:interp}, to get
  \begin{gather*}
    \abs{u-u_p}_{k,\infty,\Oreg}
       \leq \abs{u-I_p u}_{k,\infty,\Oreg} + \abs{I_p u-u_p}_{k,\infty,\Oreg}.
  \end{gather*}
  The term $\abs{u-I_p u}_{k,\infty,\Oreg}$ is handled by
  Proposition~\ref{prop:interp}, while for the second term
  an inverse inequality gives
  \begin{gather*}
    \abs{I_p u-u_p}_{k,\infty,\Oreg}
       \lesssim \abs{I_p u-u_p}_{k,2,\Oreg}
       \le \abs{I_p u-u}_{k,2,\Oreg} + \abs{u-u_p}_{k,2,\Oreg}.
  \end{gather*}
  Application of Propositions~\ref{prop:interp} and~\ref{thm:main}
  completes the proof.
\end{proof}
%%%%%%%%%%%%%%%%%%%%%%%%%%%%%%%%%%%%%%%%%%%

\begin{remark}
  Theorem \ref{thm:C1balanced} and Corollary \ref{max_norm_C1} establish the
  conjecture \eqref{C1} made in \cite{PC}.
\end{remark}

%%%%%%%%%%%%%%%%%%%%%%%%%%%%%%%%%%%%%%%%%

\section{Two-dimensional problems}
\label{sec:2D}
In two dimensional smooth domains $\Omega \subset \mathbb{R}^2$, we would have great difficulty constructing
a $C^1$	 conforming FEM; instead, a $C^0$ mixed formulation is usually preferred. Before we present the problem
under consideration along with existing results, we comment on the differences between the $C^1$ method
and the $C^0$ mixed method, beyond the obvious ones. Specifically, we want to point out that in order to use
the previous methodology, we should define two special representatives, since in the mixed formulation 
we will be seeking two unknown functions -- the second one will be auxiliary. Moreover, for the extra function we only need 
to use the $L^2$ projection, while for $u$ we (again) utilize the $H^1$ projection.

We consider the following problem\footnote{Henceforth, the symbols that appeared in the previous section will represent two-dimensional analogs.} from \cite{CFLX}: 
Find $u$ such that 
\begin{subequations}
  \label{bvp2}
  \begin{alignat}{2}
    \eps^{2}\Delta ^{2}u-b\Delta u+cu & = f & \ \ & \text{in }\Omega \subset 
  \RR^{2},  \\ % \label{de2} \\
  u=\frac{\partial u}{\partial n} & = 0  && \text{on} \ \partial \Omega ,  % \label{bc2}
  \end{alignat}
\end{subequations}
where \mbox{$0<\eps \leq 1$} is a given parameter, $\Delta $ denotes the
Laplacian (and $\Delta^{2}$ the biharmonic) operator, \mbox{$b,c>0$} are given
constants, $\Omega $ is an open bounded domain with analytic boundary
$\partial \Omega$ and $f$ is a given analytic function, which satisfies 
\begin{equation}
  \norm{\nabla ^{n}f}_{\infty ,\Omega}
     \lesssim n!\gamma_{f}^{n}\ \ \text{for} \ n=0,1,2,\dots,  \label{analytic2}
\end{equation}
for some positive constant $\gamma _{f}$ independent of $\eps$.
Here we have used the short\-hand notation 
\begin{equation*}
  \abs{\nabla ^{n}f}^{2}
     \coloneqq \sum_{\abs{\alpha}=n} 
     \frac{\abs{\alpha}!}{\alpha !}\abs{D^{\alpha} f}^{2},
     % = \sum_{\beta _{1},...,\beta_{n}=1}^{2}\left\vert D^{\beta_{1} \cdots \beta _{n}}f\right\vert ^{2},
\end{equation*}
with $D^{\alpha}f$ denoting the generalized derivative of $f$ with respect
to the multiindex $\alpha\in\NN_0^2$.
%\footnote{
%A (standard, conforming) variational formulation of~\eqref{bvp2} reads:
%Find \mbox{$u\in H_{0}^{2}\left(\Omega \right)$} such that 
%%
%\begin{gather}
%  \label{weak}
%  \eps^{2}\scal{\Delta u}{\Delta v} + b \scal{\nabla u}{\nabla v} +
%    c \scal{u}{v} = \scal{f}{v} \ \ \forall \ v\in H_{0}^{2}\left( \Omega \right).
%\end{gather}
%%
%Associated with the above problem is the \emph{energy norm} 
%%
%\begin{equation}
%\norm{u}_{E}^{2} \coloneqq \eps^{2}\scal{\Delta u}{\Delta u}
%   + b\scal{\nabla u}{\nabla u}
%   + c\scal{u}{u}.  \label{energy2}
%\end{equation}
%}

\begin{figure}[h]
  \centerline{
  \includegraphics[width=0.6 \textwidth]{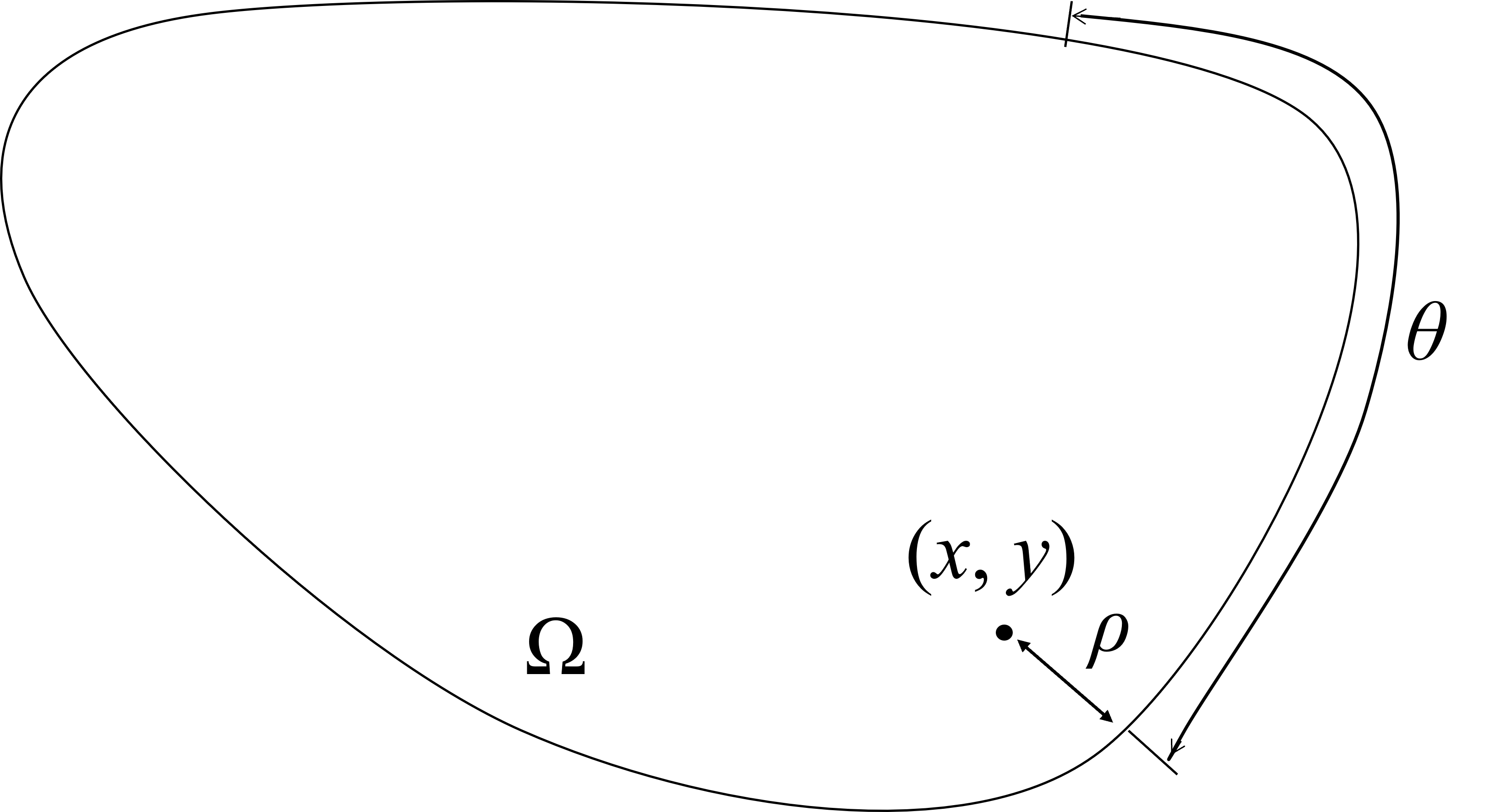}}
  \caption{Boundary fitted coordinates.}
  \label{fig0}
\end{figure}

The regularity of the solution to~\eqref{bvp2} was studied in \cite{PC, MX},
and will be described in Proposition \ref{prop:reg} below. 
To this end, define \emph{boundary fitted coordinates} \mbox{$(\rho,\theta)$}
in a neighborhood of the boundary as follows:
Let \mbox{$\bigl(X(\theta),Y(\theta)\bigr)$},
\mbox{$\theta \in [0,L]$} be a parametrization of $\partial \Omega$ 
by arclength and let $\Omega _{0}$ be a tubular
neighborhood of $\partial \Omega $ in $\Omega $.
For each point \mbox{$z=(x,y)\in\Omega_{0}$} there is a unique nearest point
\mbox{$z_{0}\in \partial \Omega$}, so with $\theta $ the arclength parameter
(with counterclockwise orientation), we set \mbox{$\rho=\abs{z-z_{0}}$} which
measures the distance from the point $z$ to $\partial \Omega$.
Explicitely, 
let $\rho _{0}>0$ be less than the minimum radius of curvature of $\partial \Omega$ and
set
\begin{gather}
  \Omega_{0} \coloneqq
     \Bigl\{ z-\rho \overrightarrow{n}_{z} \colon
       z\in \partial \Omega, \ \rho\in (0,\rho_{0}) \Bigr\},  \label{Omega0}
\end{gather}
where $\overrightarrow{n}_{z}$ is the outward unit normal
at \mbox{$z\in \partial \Omega$}, and 
\begin{equation*}
  x=X(\theta )-\rho Y'(\theta ),\ y=Y(\theta )+\rho X'(\theta),
\end{equation*}
with $\rho \in (0,\rho _{0})$, $\theta \in (0,L)$; see Figure \ref{fig0}.

\begin{proposition}
  \label{assumption}
  \begin{subequations}\label{AAA}
  The solution $u$ of~\eqref{bvp2} can be decomposed into a smooth part $u^{S}$,
  a boundary layer part $u^{BL}$ and a remainder $u^R$, viz. 
  \begin{equation}
    u = u^{S}+\chi u^{BL}+u^R,%\label{AA1}
  \end{equation}
  where $\chi $ is a smooth cut-off function, satisfying 
  \begin{equation*}
    \chi = \begin{cases}
        1 & \text{for} \ 0<\rho <\rho _{0}/3, \\ 
        0 & \text{for} \ \rho >2\rho _{0}/3.
    \end{cases}
  \end{equation*}
  Moreover, there exist constants $K_{1},K_{2},\omega ,\delta >0$, independent
  of $\eps$, such that 
  \begin{gather}
    \label{AA2}
    \norm{D^{n}u^{S}}_{0,\Omega }
       \lesssim 
          \abs{n}! \, K_{1}^{\abs{n}}\ \ \forall \ n\in \NN_{0}^{2}, \\
    \label{AA3}
    \abs{\frac{\partial ^{m+n}u^{BL}(\rho ,\theta)}{\partial \rho^{m} \partial \theta ^{n}}}
    \lesssim n!\ K_{2}^{m+n} \eps^{1-m} \E^{-\omega \rho /\eps}\ \ \forall \
        m,n\in \NN, \ (\rho,\theta) \in \overline{\Omega}_0,
    \intertext{and}
    % \label{AA4}
    \eps \abs{u^R}_{2,\Omega} + \norm{u^R}_{1,\Omega}\lesssim \E^{-\delta/\varepsilon}.
  \end{gather}
  Finally, there exist constants $C,K>0$, depending only on the data, such
  that 
  \begin{equation}
    % \label{AA6}
    \norm{D^{n}u}_{0,\Omega } \leq CK^{\abs{n}}
      \max \left\{ \abs{n}^{\abs{n}}, \eps^{1-\abs{n}}\right\} \ \
         \forall \ n\in \NN_{0}^{2}.
  \end{equation}
  \end{subequations}
\end{proposition}
\begin{proof}
  All estimates were shown in \cite{PC}, except for \eqref{AA2} which
  appears in \cite{MX2}.
\end{proof}

Letting\footnote{The fact that $w \in H^2$ follows from the smoothness of $\pt\Omega$.}
\mbox{$w\coloneqq\eps \Delta u\in H^{2}(\Omega)$},
we shall study the \emph{mixed} formulation of finding 
\mbox{$(u,w)\in H_{0}^{1}(\Omega )\times H^{1}(\Omega)$}
such that 
\begin{subequations}
  \label{mixed}
  \begin{alignat}{2}
    \eps \scal{\nabla u}{\nabla \phi} + \scal{w}{\phi} & = 0 & \ \ &
        \forall \ \phi \in H^{1}(\Omega), \\ 
    b\scal{u}{\nabla} + c\scal{u}{\psi} - \eps \scal{\nabla w}{\nabla}
    & = \scal{f}{\psi} && \forall \ \psi \in H_{0}^{1}(\Omega ).
  \end{alignat}
\end{subequations}
We define the bilinear form
\begin{gather}
  \label{Buv2}
    \m{A}_{\eps} \bigl((u,w),(\psi ,\phi )\bigr)
     \coloneqq
        \eps \scal{\nabla u}{\nabla \phi} +
        \scal{w}{\phi}
    + b \scal{\nabla u}{\nabla \psi} +
        c\scal{u}{\psi} - \eps \scal{\nabla w}{\nabla \psi}.
\end{gather}
Then  the \emph{energy norm} of $u$ is given by
\begin{gather*}
  % \label{norm-coercive}
  \norm{u}_{E}^{2} =
     \enorm{(u,\eps \Delta u)}^{2}
       = \eps^{2} \norm{\Delta u}_{0,\Omega }^{2}
          + b\norm{\nabla u}_{0,\Omega }^{2}
          + c\norm{u}_{0,\Omega}^{2} =  \enorm{(u,w)}^{2}.
\end{gather*}

A regularity result for $w$, analogous to Proposition \ref{assumption}, is given below.
\begin{corollary}
  \begin{subequations}
  \label{B0}
  Let \mbox{$w = \eps \Delta u$} and assume Proposition~\ref{assumption} holds.
  Then
  \begin{equation}
    w=w^{S}+\chi w^{BL}+w^R,  \label{B1}
  \end{equation}
  where $\chi $ is the same cut-off function as in Proposition \ref{assumption}.
  Moreover, there exist constants $K_{3},K_{4},\beta ,\gamma >0$, independent
  of $\eps $, such that 
  \begin{gather}
    % \label{B2}
    \norm{D^{n}w^{S}}_{0,\Omega } \lesssim \abs{n}!\ K_{3}^{\abs{n}}
      \ \ \forall \ n\in \NN_{0}^{2}, \\
    \label{B3}
    \abs{\frac{\pt^{m+n}w^{BL}(\rho,\theta )}{\pt\rho^{m} \pt\theta ^{n}}}
      \lesssim n!K_{4}^{m+n}\eps^{-m}
      \E^{-\beta\rho/\eps}\ \ \forall \ m,n\in\NN, \ (\rho,\theta) \in \overline{\Omega}_0,
    \intertext{and}
    % \label{B4}
     \norm{w^R}_{1,\Omega} \lesssim \E^{-\gamma/\eps}.
  \end{gather}
  Finally, there exists a constant $K>0$, depending only on the data, such that 
  \begin{gather}
    % \label{B5}
    \norm{D^{n}w}_{0,\Omega }\lesssim K^{\abs{n}} 
       \max \left\{ \abs{n}^{\abs{n}}, \eps^{-\abs{n}}\right\} \ \ \forall \ n\in \NN_{0}^{2}.
  \end{gather}
  \end{subequations}
\end{corollary}
\begin{proof}
  Follows from Proposition \ref{assumption} and the definition
  of $w=\eps \Delta u$.
\end{proof}

We close this section with the following result which will be used in the
sequel.
\begin{lemma}\label{BLsize}
  Let $u^{BL}$ and $w^{BL}$ satisfy~\eqref{AA3} and \eqref{B3}, respectively.
  With $\Omega_0$ given by~\eqref{Omega0}, we have
  \begin{gather*}
    \eps^{-1} \norm{u^{BL}}_{0,\Omega \setminus \Omega_0}
       +
    \norm{w^{BL}}_{0,\Omega \setminus  \Omega_0} \lesssim \eps^{1/2} \E^{-\sigma/ \eps},
  \end{gather*}
  for some positive constant $\sigma$, independent of $\eps$.
\end{lemma}
\begin{proof}
  This follows from direct calculations.
\end{proof}

\subsection{Discretization by a mixed $rp$-FEM}% \label{mesh}}
In order to define our finite dimensional discretization spaces
\mbox{$V_{1}^{N}\subset H_0^1(\Omega)$} and \mbox{$V_{2}^{N}\subset H^1(\Omega)$},
we let 
\mbox{$\Delta \coloneqq \left\{\Omega _{i}\right\} _{i=1}^{N}$} be a mesh
consisting of curvilinear quadrilaterals, subject to standard conditions
(see, e.g. \cite{melenk}) and associate with each $\Omega _{i}$ a bijective
mapping \mbox{$M_{i}\colon S_{ST} \rightarrow \overline{\Omega }_{i}$},
where \mbox{$S_{ST}=[0,1]^{2}$} denotes the reference square.
With \mbox{$\m{Q}_{p}(S_{ST})$} the space of polynomials of degree $p$ (in each
variable) on $S_{ST}$, we define 
\begin{align*}
  \m{S}^{p}(\Delta )
    & \coloneqq \Bigl\{ u\in H^{1}\left( \Omega \right) \colon
           u\vert_{\Omega_{i}}\circ M_{i}\in \m{Q}_{p}(S_{ST}), \ \ i=1,\dots,N\Bigr\}, \\
  \m{S}_{0}^{p}(\Delta)
    & \coloneqq \m{S}^{p}(\Delta )\cap H_{0}^{1}(\Omega ).
\end{align*}
The mesh $\Delta$ is chosen following the construction in \cite{CFLX,MX}.
We denote by $\Delta _{A}$ a \emph{fixed} (asymptotic) mesh
consisting of curvilinear quadrilateral elements $\Omega _{i}$,
\mbox{$i=1,\dots,N_{1}$}.
These elements $\Omega _{i}$ are the images of the reference square $S_{ST}$
under the element mappings $M_{A,i}$, \mbox{$i=1,\dots,N_{1}\in \NN$}. 
Moreover, the element mappings $M_{A,i}$ are assumed to be analytic (with
analytic inverses). 
We also assume that the elements do not have a single vertex on the
boundary~$\pt\Omega$, but only complete, single edges.
For convenience, we number the elements along the boundary first,
i.e., $\Omega_{i}$, \mbox{$i=1,\dots,N_{2}<N_{1}$} for some \mbox{$N_{2}\in \NN$}.

We now give the definition of the layer-adapted mesh from \cite{MX}, almost
verbatim.
\begin{definition}%[Spectral Boundary Layer Mesh $\Delta _{BL}$]
  \label{SBL-2D}
  Given parameters \mbox{$\kappa>0$}, \mbox{$p\in\NN$}, \mbox{$\eps\in (0,1]$}
  and the (asymptotic) mesh $\Delta _{A}$, the \emph{Spectral Boundary Layer
  Mesh} $\Delta _{BL}$ is defined as follows:
  \begin{enumerate}
    \item
      If \mbox{$\kappa p\varepsilon \geq 1/2$} then we are in the asymptotic
      range of p and we use the mesh $\Delta _{A}$.

   \item
     If \mbox{$\kappa p\varepsilon <1/2$}, we need to define so-called needle
     elements.
     We do so by splitting the elements $\Omega_{i}$, \mbox{$i=1,\dots,N_{2}$}
     into two elements $\Omega _{i}^{need}$ and $\Omega_{i}^{reg}$.
     To this end, split the reference square $S_{ST}$ into two rectangles 
     \begin{equation*}
       S^{need} = \left[0,\kappa p\eps \right] \times [0,1],\quad
       S^{reg}  = \left[\kappa p\eps, 1\right] \times [0,1],
     \end{equation*}
     and define the elements $\Omega _{i}^{need}$ and $\Omega _{i}^{reg}$ as
     the images of these two rectangles under the element map $M_{A,i}$ and the
     corresponding element maps as the concatination of the affine maps 
     \begin{alignat*}{5}
       & A^{need} && \colon S_{ST} \rightarrow S^{need},
                && \quad (\xi ,\eta )\rightarrow (\kappa p\varepsilon \xi ,\eta ), \\
       & A^{reg}  && \colon S_{ST} \rightarrow S^{reg},
                && \quad (\xi ,\eta )\rightarrow (\kappa p\varepsilon +
     (1-\kappa p\varepsilon )\xi ,\eta )
     \end{alignat*}
     with the element map $M_{A,i}$, i.e., \mbox{$M_{i}^{need}=M_{A,i}\circ A^{need}$}
     and $M_{i}^{reg}=M_{A,i}\circ A^{reg}$.
  \end{enumerate}
\end{definition}

\begin{figure}[h]
\begin{center}
\includegraphics[width=0.3 \textwidth]{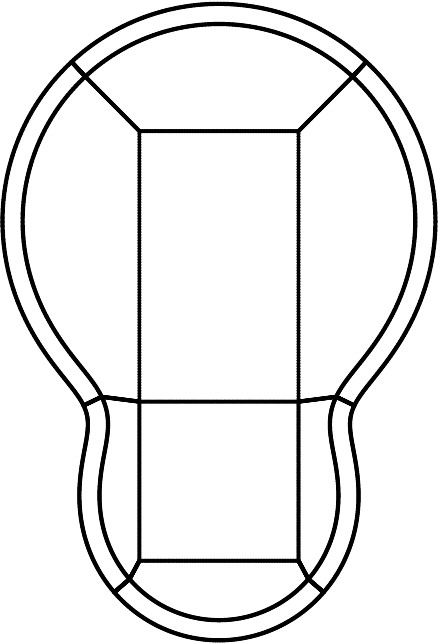}
\end{center}
\caption{Example of the SBL mesh.}
\label{fig1}
\end{figure}

In total, the mesh $\Delta_{BL}$ consists of \mbox{$N=N_{1}+N_{2}$}
elements if \mbox{$\kappa p\varepsilon <1/2$}.
By construction, the resulting mesh 
\begin{gather*}
  \Delta _{BL}
    = \Bigl\{\Omega_{1}^{need},\dots,
             \Omega_{N_{1}}^{need},\Omega _{1}^{reg},\dots,
             \Omega _{N_{1}}^{reg},
             \Omega_{N_{1}+1},\dots,\Omega _{N}\Bigr\},
\end{gather*}
is a regular admissible mesh.
In Figure \ref{fig1}, we show an example of such a mesh from \cite{CFLX}.

Next, we take \mbox{$V_{1}^{N}=\mathcal{S}_{0}^{p}(\Delta_{BL})$} and
\mbox{$V_{2}^{N}=\mathcal{S}^{p}(\Delta_{BL})$},
and seek an approximation $\left(u_{p,\Delta},w_{p,\Delta}\right)\in
V_1^N\times V_2^N$ of~\eqref{mixed}
such that
\begin{gather}
  \label{discrete2}
  \m{A}_{\eps} \left(\left(u_{p,\Delta},w_{p,\Delta}\right)),(\psi ,\phi )\right)
    = \scal{f}{\psi} \ \ \forall \ (\phi,\psi)\in V_1^N\times V_2^N.
\end{gather}
Subtracting~\eqref{discrete2} from~\eqref{mixed}, we get 
\begin{gather*}
  % \label{orthogonal}
  \m{A}_{\eps} \left( (u-u_{p,\Delta},w-w_{p,\Delta}),(\psi ,\phi )\right) =0
     \ \ \forall \ (\psi ,\phi)\in V_{1}^{N}\times V_{2}^{N}. 
\end{gather*}

We have the following interpolation estimates, which were based on the Gauss-Lobbato interpolant from \cite{melenk}.
\begin{lemma}[\cite{CFLX}]% \label{lem:CFLX}
  Let $\left(u,w\right)$ be the solution to~\eqref{mixed} and assume
  that~\eqref{analytic2} holds.
  Then there exist constants \mbox{$\kappa_{0},\kappa _{1}, \beta >0$}
  independent of both \mbox{$\eps \in (0,1]$} and of \mbox{$p\in \NN$},
  such that the following is true: For every $p$ and every
  \mbox{$\kappa \in (0,\kappa _{0}]$} with \mbox{$\kappa p\geq \kappa _{1}$},
  there exist \mbox{$I _{p}u\in \mathcal{S}_{0}^{{p}}(\Delta _{BL})$} and 
  \mbox{$J_{p}w\in \mathcal{S}^{{p}}(\Delta _{BL})$} such that 
  \begin{gather*}
    \norm{u-I_{p}u}_{L^{\infty}(\Omega)}
      + \norm{\nabla (u-I_{p}u)}_{L^{\infty}(\Omega)}
      + \norm{w-J_{p}w}_{L^{\infty}(\Omega)}
      + \eps^{1/2} \norm{\nabla(w-J_{p}w)}_{0 ,\Omega }
          \lesssim \E^{-\beta p}.
  \end{gather*}
\end{lemma}
With
\begin{gather*}
  \Omega_{BL} = \bigcup_i \Omega_i^{need} \, ,
      \quad \Oreg = \bigcup_i \Omega_i^{reg}
\end{gather*}
and \mbox{$\abs{\Omega_{BL}} = \ord{\kappa p \eps}$},
\mbox{$\abs{\Oreg} = \ord{1}$}, we see that additionally there holds
\begin{equation}
  \label{BLregion_bound}
  \norm{\nabla (u-I_{p} u)}_{0 ,\Omega_{BL}}
    + %\lesssim \varepsilon^{1/2} e^{-\beta p } \; , \;
  \norm{w-J_{p} w}_{0 ,\Omega_{BL}} \lesssim \eps^{1/2} \E^{-\beta p }
\end{equation}

The following was the main result of \cite{CFLX}.
\begin{theorem}[\cite{CFLX}]
  \label{theorem_main}
  % Let $(u,w)\in H_{0}^{1}(\Omega )\times H^{1}(\Omega ),$ 
  % $\left( u_{FEM},w_{FEM}\right) \in \mathcal{S}^p_0(\Delta_{BL}) \times \mathcal{S}^p(\Delta_{BL})$ 
  Let $(u,w)$ and $\left( u_{p,\Delta},w_{p,\Delta}\right)$ be the solutions to
  \eqref{mixed} and \eqref{discrete2}, respectively.
  Assume Proposition \ref{assumption} holds. 
  Then there exists a positive constant $\beta$, independent of $\eps$, such that
  \begin{gather*}
    \enorm{\left( u-u_{p,\Delta},w-w_{p,\Delta}\right)} \lesssim \E^{-\beta p}.
  \end{gather*}
\end{theorem}

Note that the norm used in Theorem~\ref{theorem_main} is not correctly
balanced, because
\begin{equation*}
  \enorm{(u^{BL},\eps \Delta u^{BL})} = \ord{\eps^{1/2}}\,,
     \ \ \text{but} \ \ \enorm{(u^{S},\eps \Delta u^{S})} = \ord{1}\,,
   \ \ (\eps \to 0).
\end{equation*}
Just like in the one-dimensional case,
as \mbox{$\eps \rightarrow 0$}, ``the energy norm does
not see the layers'' \cite{FR, FR3}. This is due to
the fact that the weight on $\norm{w}_{0,\Omega}$ is not the
appropriate one.
A more suitable, correctly balanced norm is 
\begin{gather}
  \label{balanced}
  \enorm{(u,w)}_{B}^{2}
     \coloneqq \eps^{-1} \norm{w}_{0,\Omega}^{2}
       + b \norm{\nabla u}_{0,\Omega }^{2}
       + c \norm{u}_{0,\Omega }^{2},
\end{gather}%
since then
\begin{equation*}
  \enorm{\left(u^{BL},\eps \Delta u^{BL}\right)}_B
     = \ord{1} = 
  \enorm{\left(u^{S},\eps \Delta u^{S}\right)}_B, \ \ \ (\eps\to0).
\end{equation*}
Numerical experiments using the above norm have been reported in \cite{CFLX},
suggesting that indeed the method yields uniform exponential convergence in the
balanced norm $\enorm{\cdot}_B$.
Unfortunately, the bilinear form $\m{A}_{\eps}$ given by~\eqref{Buv2}, is not
uniformly coercive with respect to the norm~\eqref{balanced}.
To circumvent this obstacle, we will use the methodology derived in
Section~\ref{sec:1D}.

%%%%%%%%%%%%%%%%%%%%%%%%%%%%%%%%%%%%%%%%%%%%%
%----------------------
\subsection{Balanced norm estimates}
\label{estimates}

Recall the decomposition~\eqref{B0}, set \mbox{$w^{SR}\coloneqq w^S+w^R$}
and define the $L^2$ projection of $w^{SR}$ on
the regular region $\Oreg$, denoted by
\mbox{$\pi^1_p w^{SR} \in V^N_2|_{\Oreg}$}, via
\begin{gather*} % \label{pi_p2D}
  \scal{\pi^1_p w^{SR}}{v}_{\Oreg}=0 \quad \forall \; v 
    \in V^N_2|_{\Oreg}.
\end{gather*}
The projection is uniquely defined and consequently, we have the following:
\begin{gather}
  \label{pi_p2}
   \norm{\pi^1_p w^{SR}-w^{S}}_{0, \Omega_{REG}}
     + \norm{\pi^1_p w^{SR}-w^{SR}}_{L^{\infty}(\Omega_{REG})}
          \lesssim \E^{-\kappa p}
\end{gather}
where the definition of the projection $\pi^1_p$ and Lemma~\ref{BLsize}
were used, along with $\kappa p \eps < 1/2$, which is when we use the SBL mesh.

Let analogously \mbox{$u^{SR}\coloneqq u^S+u^R$} and define
\mbox{$\pi_p^2 u^{SR} \in V^N_2|_{\Omega_{REG}}$} as the
 weighted $H^1$ projection
of $u^{SR}$ on  $\Omega_{REG}$:
\begin{align*} %\label{pi_p2_2D}
  b \scal{\nabla ( \pi^2_p  u^{SR}-u^{SR})}{v}_{\Omega_{REG}}
    + c \scal{\pi^2_p u^{SR}-u^{SR}}{v}_{\Omega_{REG}}=0 \quad \forall \; v 
  \in V^N_2|_{\Omega_{REG}},
\end{align*}
and we have
\begin{gather*}% \label{lem:pi_p3}
  \norm{\pi^2_p u^{SR}-u^{SR}}_{1, \Omega_{REG}}
    + \norm{\pi^2_p \nabla (u^{SR}-u^{SR})}_{L^{\infty}(\Omega_{REG})}
      \lesssim \E^{-\beta p}.  %\label{lem:pi_p4}
\end{gather*}

We next define an auxiliary function $\chi_2 \in  \mathcal{P}_1$ by the conditions
\begin{gather*}
  \chi_2(0)=0, \quad \chi_2(\kappa p \eps)=1.
\end{gather*}
A direct calculation establishes the following bounds on $\chi_2$:
\begin{gather}\label{bound-chi2D}
  \norm{\chi_2}_{0,\Omega_{BL}} \lesssim (\kappa p \varepsilon)^{1/2}, \quad
  \norm{\chi_2'}_{0,\Omega_{BL}} \lesssim (\kappa p \varepsilon)^{-1/2}. \quad
\end{gather}

We are now in a position to define \emph{special representatives} $\tilde{w}\in \m{S}^p$
of $w$ and $\tilde{u}\in \m{S}_0^p$ of $u$.
Recalling the decomposition~\eqref{B1}, define

\begin{subequations}
\begin{align}\label{wt}
  \tilde{w}(x) & \coloneqq
    \begin{cases}J_p w-  \chi_2(\rho) \left(w^{SR}+w^{BL}-\pi^1_p w^{SR}\right)|_{\partial \Omega_{REG}} & \text{ in } \Omega_{BL}, 
\\ \pi_p^1 w^{SR} & \text{ in } \Omega_{REG}.  
  \end{cases} \\
  \label{ut}
  \tilde{u}(x) & \coloneqq
    \begin{cases}I_p u-  \chi_2(\rho) \left(u^{SR}+u^{BL}-\pi^2_p u^{SR}\right)|_{\partial \Omega_{REG}} & \text{ in } \Omega_{BL}, 
\\ \pi_p^2 u^{SR} & \text{ in } \Omega_{REG}.  
\end{cases}
\end{align}
\end{subequations}
We have the following.

\begin{lemma}\label{lem:u-ut2D}
  There exist a positive contant $\beta$ such that the following bounds holds true:
  \begin{subequations}
  \begin{align} \label{w-wt2D:bl}
    \norm{w - \tilde{w}}_{0,\Omega_{BL}} & \lesssim  \eps^{1/2} \E^{-\beta p}, \\
  \label{w-wt2D:reg}
    \norm{w - \tilde{w}}_{0,\Omega_{REG}} & \lesssim  \E^{-\beta p}, \\
     \label{u-ut2D:bl}
     \norm{u - \tilde{u}}_{0,\Omega_{BL}} & \lesssim  \eps^{1/2} \E^{-\beta p}, \\
  \label{u-ut2D:reg}
    \norm{u - \tilde{u}}_{0,\Omega_{REG}} & \lesssim  \E^{-\beta p}.
  \end{align}
  \end{subequations}
\end{lemma}
\begin{proof} We study the two regions $\Olay$ and $\Oreg$ separately. 

  \paragraph{$\Olay$:}
  We have for $w-\tilde{w}$,
  \begin{gather*}
    \norm{w - \tilde{w}}_{0,\Omega_{BL}}
      \le \norm{w- J_p w}_{0,\Omega_{BL}}
      + \norm{\chi_2}_{0,\Omega_{BL}} \left\{ \norm{ w^{BL}}_{L^{\infty}(\partial \Oreg)}
              + \norm{\pi^1_p w^{SR} - w^{SR}}_{L^{\infty}(\Omega_{BL})}\right\}\,.
  \end{gather*}
  Equation \eqref{BLregion_bound} takes care of the first term above,
  and equations~\eqref{AA3},~\eqref{pi_p2} and~\eqref{bound-chi2D},
  allow us to bound the remaining terms, as desired.
  This shows \eqref{w-wt2D:bl}.
  Similarly, for $u-\tilde{u}$ we have
  \begin{gather*}
    \norm{u - \tilde{u}}_{0,\Omega_{BL}}
      \le \norm{u- I_p u}_{0,\Omega_{BL}}
      + \norm{\chi_2}_{0,\Omega_{BL}} \left\{ \Vert{ u^{BL} \Vert_{L^{\infty}(\partial \Oreg)} }
              + \norm{\pi^2_p u^{SR} - u^{SR}}_{L^{\infty}(\Omega_{BL})}\right\}.
  \end{gather*}
  Eq. \eqref{u-ut2D:bl} follows from the same equations.
  
  \paragraph{$\Oreg$:}
  We have $w=w^{BL} + w^{SR}$.
  Therefore,
  \begin{gather*}
    \norm{w - \tilde{w}}_{0,\Oreg} \le
      % \norm{u_\mathrm{bl}}_{\Oreg} + \norm{u_\mathrm{reg} - \tilde{u}}_{\Oreg} =
      \norm{w^{BL}}_{0,\Oreg}
        + \norm{w^{SR} - \pi^1_p w^{SR} }_{0,\Oreg} .
  \end{gather*}
  Lemma \ref{BLsize} and \eqref{pi_p2} then give \eqref{w-wt2D:reg}. Eq. \eqref{u-ut2D:reg} follows in
  the same fashion.
\end{proof}

We now have the necessary tools for proving the following.

\begin{theorem}% \label{thm:balanced}
  Let $(u,w)$ be the solution of \eqref{mixed} and let $(u_{p,\Delta},w_{p,\Delta})$
  be the solution of \eqref{discrete2}, based on the Spectral Boundary Layer mesh
  of Definition \ref{SBL-2D}.
  Then, there exists a positive constant $\beta$, independent of $\eps$, such that
  \begin{gather*}
    \enorm{(u - u_{p,\Delta}, w - w_{p,\Delta})}_{B} \lesssim \E^{-\beta p}.
  \end{gather*}
\end{theorem}
\begin{proof}
  Let $\tilde{w}, \tilde{u}$ be defined by \eqref{wt} and \eqref{ut}, respectively. 
  We have from \eqref{balanced} and the triangle inequality,
  \begin{align*}
    \enorm{(u-u_{p,\Delta},w-w_{p,\Delta})} _{B}
    & \leq \enorm{(u-\tilde{u},w-\tilde{w})} _{B}
           + \enorm{(\tilde{u}-u_{p,\Delta},\tilde{w}-w_{p,\Delta})} _{B} \\
    & \lesssim \eps^{-1/2} \norm{w - \tilde{w}}_0
          + \norm{u - \tilde{u}}_1 + \eps^{-1/2} \norm{\tilde{w}-w_{p,\Delta}}_0
          + \norm{\tilde{u} - u_{p,\Delta}}_1 .
  \end{align*}
  We only need to treat the terms $\eps^{-1/2} \norm{\tilde{w}-w_{p,\Delta}}_0$
  and $\norm{\tilde{u}-u_{p,\Delta}}_1$ since the rest can be handled by
  Theorem \ref{theorem_main} and Lemma \ref{lem:u-ut2D}. 

  Let $\eta \coloneqq w_{p,\Delta} - \tilde{w}$ and
  $\xi \coloneqq u_{p,\Delta} - \tilde{u}$.
  We have by Galerkin orthogonality,
  \begin{gather}\label{GO}
    \enorm{(\eta, \xi)} =
       \scal{w - \tilde{w}}{\eta} + b \scal{\nabla (u - \tilde{u})}{\nabla \xi}
         +  c \scal{ u - \tilde{u}}{\xi}.
  \end{gather}

  We will consider the two regions separately.

  \paragraph{$\Omega_{BL}$:} We first note that
  \begin{gather*}
    \norm{\chi_2}_{0,\Omega_{\mathrm{BL}}}
      \left\{\norm{w^{BL}}_{L_\infty(\pt\Omega_{REG})} +
        \norm{\pi^1_p w^{SR}-w^{SR}}_{L_\infty(\pt\Omega_{REG})} \right\}
      \lesssim \eps^{1/ 2} \E^{-\kappa p}.
  \end{gather*}
  Then, we combine the above with \eqref{BLregion_bound},
   and use the Cauchy-Schwarz inequality, along with 
  the defintion of $\hat{w}$, to obtain
  \begin{gather*}
    \scal{w-\hat{w}}{\eta}_{\Omega_{BL}}
       \lesssim \eps^{1/ 2} \E^{-\beta p} \norm{\eta}_{\Omega_{BL}}.
  \end{gather*}
  Similarly, we first note that
  \begin{gather*}
    \abs{\chi_2'}_{0,\Omega_{\mathrm{BL}}}
      \left\{\norm{u^{BL}}_{L_\infty(\pt\Omega_{REG}} +
        \norm{\pi^2_p u^{SR}-u^{SR}}_{L_\infty(\pt\Omega_{REG}} \right\}
           \lesssim \eps^{1/2} \E^{-\kappa p},
  \end{gather*}
  hence, as above, we have
  \begin{gather*}
    b\scal{\nabla (u - \tilde{u})}{\nabla \xi}_{\Omega_{BL}}
       +c\scal{u - \tilde{u}}{\xi}_{\Omega_{BL}}
      \lesssim \eps^{1/2} \E^{-\beta p} \norm{\xi}_{1,\Omega_{BL}}.
  \end{gather*}

  \paragraph{$\Omega_{REG}$:} We have by the definition of $\hat{w}$,
  \begin{align*}
    \scal{w-\hat{w}}{\eta}_{\Omega_{REG}}
    = \scal{w_{BL} + w^{SR} - \pi_p w^{SR}}{\eta}_{\Oreg}
      = \scal{w^{BL}}{\eta}_{\Omega_{REG}}
        \lesssim \eps^{1/2} \E^{-\kappa p} \norm{\eta}_{0, \Oreg}.
  \end{align*}
  Similarly,
  \begin{align*}
    & b \scal{\nabla (u - \tilde{u})}{\nabla \xi}_{\Oreg}
      + c \scal{u - \tilde{u}}{\xi}_{\Oreg} \\
    & \qquad = b \scal{\nabla u^{BL}}{\nabla \xi}_{\Oreg}
      + c \scal{u^{BL}}{\xi}_{\Oreg}
      \lesssim \norm{u^{BL}}_{1,\Olay} \norm{\xi}_{1,\Olay}
      \lesssim \eps^{1/2} \E^{-\beta p} \norm{\xi}_{1,\Olay}\,.
  \end{align*}

  Combining the above with \eqref{GO}, we see that
  \begin{gather*}
    \enorm{(\eta, \xi)}^2 \lesssim \eps^{1/2} \E^{-\beta p} \enorm{(\eta, \xi)}\,,
  \end{gather*}
  from which we have
  \begin{gather*}
    \eps^{-1/2} \norm{w_{p,\Delta} - \tilde{w}}_{0,\Omega}
      + \norm{u_{p,\Delta} - \tilde{u}}_{1,\Omega} \lesssim \E^{-\beta p}.
  \end{gather*}
  This completes the proof.
\end{proof}
%

%\subsection{Maximum norm estimates}
%
%Since we have the balanced norm estimate, we may obtain a maximum norm
%estimate in the following.
%
%\begin{corollary}
%  Let the assumptions of Theorem \ref{thm:balanced} be true.
%  Then there exists a positive constant $\beta$ independent of $\eps$,
%  such that
%  %
%  \begin{gather*}
%    \norm{w - w_{p,\Delta}}_{L^{\infty}(\Omega)}
%     + \norm{u - u_{p,\Delta}}_{L^{\infty}(\Omega)}
%       \lesssim \E^{-\beta p},
%  \end{gather*}
%\end{corollary}
%
%\begin{proof}
%  We have from the triangle inequality,
  %
%  \begin{gather*}
%    \norm{w - w_{p,\Delta}}_{L^{\infty}(\Omega)}
%      \le \norm{w - J_p w}_{L^{\infty}(\Omega)}
%            + \norm{J_p w - w_{p,\Delta}}_{L^{\infty}(\Omega)},
%  \end{gather*}
%  %
%  with $J_p$ the interpolant of Lemma \ref{lem:CFLX}.
%  So it suffices to bound the term $\norm{J_p w - w_{p,\Delta}}_{L^{\infty}(\Omega)}$.
%  In addition, we only need to do so in the layer region, since in $\Omega_{REG}$
%  standard inverse estimates give the desired result. To that end, consider
%  %
%  \begin{gather*}
%    \norm{J_p w - w_{p,\Delta}}_{L^{\infty}(\Omega^{need}_i)}
%      \lesssim \frac{p^2}{\sqrt{\kappa p \eps}} \norm{J_p w - w_{p,\Delta}}_{0,\Omega^{need}_i}
%  \end{gather*}
%  %
%  where we used the polynomial inverse estimate of \cite[Thm. 4.76]{schwab}.
%  Equation \eqref{BLregion_bound} gives the result.
%\footnote{I doubt. Eq.~\eqref{BLregion_bound} gives interpolation error
%  bounds for $w-J_pw$, not for $J_p w - w_{p,\Delta}$}
%  The same technique works for $u$ as well.
%\end{proof}
%
\section{Conclusions}\label{concl}

We presented balanced norm estimates for singularly perturbed fourth order problems in one-dimension
as well as in smooth (analytic) domains in $\mathbb{R}^2$. As a corollary, we also got maximum norm
estimates in the one-dimensional case. The methodology is not restricted to conforming methods, and in fact we are in the process of studying balanced norm estimates for non-conforming $rp$-version methods.

Finally, we mention that this contribution proves the conjectures made in \cite{PC}, \cite{CFLX} based on
numerical evidence.

%%%%%%%%%%%%%%%%%%%%%%%%%%%%%%%%%%%%%%%%%%%

%---------------------------------------

\appendix

\section{Auxiliary results}

%The following inverse inequality is used throughout the article.
%
\begin{lemma}[Inverse inequality]
  \label{lem:inverse-1D}
  For any polynomial $q\in\mP_p$, $p\in\NN$, and any domain $D$ there holds
  \begin{gather*}
    \abs{q}_{k,D} \le \left(\frac{p!}{(p-k)!}\right)^2 \abs{D}^{-k} \norm{q}_{0,D}\,,
       k=0,1,\dots, p.
  \end{gather*}
\end{lemma}
\begin{proof}
Theorem 4.76 in \cite{schwab} gives the cases $k=0$ and $1$. Iterating on that result, we obtain
the desired inequality.

\end{proof}

\section{Proof of \eqref{interp_uS} and \eqref{interp_lay}}
We begin with the following results:
\begin{lemma}
\label{lem_aux1}
Let \mbox{$v \in C^{\infty}(I_{ST})$}, \mbox{$I_{ST}=(-1,1)$} satisfy
\begin{gather*}
  \norm{v^{(n)}}_{{\infty},I_{ST}}
    \le C_v (\gamma_v h)^n K^{n-1}, \  \forall \; n \in \mathbb{N},
\end{gather*}
for some constants \mbox{$C_v, \gamma_v>0, K\ge 1$}, and for \mbox{$h \in (0,1]$}.
Then there exists an approximation \mbox{$\mathcal{I}_p v \in \mP_p$} with
\mbox{$\mathcal{I}^{(k)}_p v(\pm 1) = v^{(k)}(\pm 1)$}, \mbox{$k=0,1$} such that
under the condition 
\begin{gather*}
  \frac{h K}{p} \leq \eta,
\end{gather*}
for some \mbox{$\eta > 0$}, there holds
\begin{gather*}
  \norm{v - \mathcal{I}_p v}_{2,I_{ST}} \lesssim  C_v K^{1/2} h \E^{-\beta p}.
\end{gather*}
\end{lemma}
\begin{proof}
  There exists (see, e.g.~\cite{PZMX}) \mbox{$\mathcal{I}_p v \in \mP_p$} with
  \mbox{$\left(\mathcal{I}_p v\right)^{(k)} (\pm 1) = v^{(k)}(\pm 1)$}, \mbox{$k=0,1$}
  and 
  \begin{gather*}
    \norm{v - \mathcal{I}_p v}^2_{2,I_{ST}}
       \leq \frac{(p-s)!}{(p+s)!} \norm{v^{(s+1)}}^2_{0,I_{ST}}, \ \ \forall \  s \in (0,p).
  \end{gather*}
  Choose $s=\lambda p$, with \mbox{$\lambda \in (0,1]$} to be selected shortly.
  Then we calculate 
  \begin{gather*}
     \frac{(p-\lambda p)!}{(p+\lambda p)!} \norm{v^{(\lambda p+1)}}^2_{0,I_{ST}}
        \le \left[ \frac{(1-\lambda)^{1-\lambda}}{ (1+\lambda)^{1+\lambda}} \right]^p
                 p^{-2\lambda p} \E^{2\lambda p + 1} C^2_v (h \gamma_v)^{2 (\lambda p+1)} K^{2\lambda p+1}
        \le \E C_v^2 h^2 \gamma_v^2 K \left( \frac{\E \gamma_v  h K}{p} \right)^{2\lambda p} q^p ,
  \end{gather*}
  with \mbox{$q = (1-\lambda)^{1-\lambda} / (1+\lambda)^{1+\lambda} \in (0,1)$}.
  Using the hypothesis $hK/p \leq \eta$ and selecting \mbox{$\lambda \in (0,1]$},
  \mbox{$\eta>0$} such that $\left( \gamma_v \E \eta \right)^{2\lambda p}\leq 1$,
  gives
  \begin{gather*}
    \frac{(p-\lambda p)!}{(p+\lambda p)!} \norm{v^{(\lambda p+1)}}^2_{0,I_{ST}}
    \lesssim C_v^2 h^2 K \gamma_v q^p,
  \end{gather*}
  from which the desired result follows.
\end{proof}

By means of affine transformations the interpolation operator $\m{I}_p$
of Lemma~\ref{lem_aux1} can be applied to all subintervals of our mesh $\Delta$
to give a global piecewise polynomial interpolant.
To this end introduce the affine transformations
\begin{gather*}
  Q_j \colon I_{ST} \to I_j
      \colon t \mapsto Q_{j}(t)
         \coloneqq\frac{1-t}{2} x_{j-1}+\frac{1+t}{2} x_j,
           \ \ j=1,\dots,N,
\end{gather*}
and define
\begin{gather}\label{interp:pw-poly}
  \left(\m{I}_{p,\Delta}v\right)\big|_{I_j}
      \coloneqq \m{I}_p\left(v\circ Q_j^{-1}\right),
           \ \ j=1,\dots,N,
\end{gather}
As there is no confusion, we drop the $\Delta$ from the notation.

\begin{corollary}
  % \label{coraux} not referenced!
  Let $j\in\{1,\dots,N\}$ and assume \mbox{$v \in C^{\infty}(\bar{I}_j)$} satisfies 
  \begin{gather*}
    \norm{v^{(n)}}_{{\infty},I} \le C_v (h_j \gamma_v)^n K^{n-1}, \ \ \forall \  n \in \mathbb{N},
  \end{gather*}
  for some constants \mbox{$C_v, \gamma_v>0$}, \mbox{$K\ge 1$}.
  Then there exist constants \mbox{$\eta, \beta > 0$} depending only on
  $\gamma_v$ such that under the condition
  \begin{gather*}
    \frac{h_j K}{p} \leq \eta,
  \end{gather*}
  the polynomial approximation $\mathcal{I}_p v$ defined in~\eqref{interp:pw-poly}
  satisfies
  \begin{gather*}
    \norm{\left(v-\mathcal{I}_p v \right)^{(k)}}_{0,I_j}
       \lesssim C_v K^{1/2} h_j^{3/2-k} \E^{-\beta p}, \ \ k=0,1,2.
  \end{gather*}
\end{corollary}
\begin{proof}
  Set \mbox{$\hat{v} = v \circ Q_j$}.
  Then
  \begin{gather*}
    \norm{\hat{v}^{(n)}}_{0,I_{ST}} \le C_v (h_j/2)^n \gamma^n_v K^{n-1},
  \end{gather*}
  and by Lemma \ref{lem_aux1} we have
  \begin{gather*}
    \norm{\hat{v} - \mathcal{I}_p \hat{v}}_{2,I_{ST}}
      \lesssim C_v K^{1/2} h_j \E^{-\beta p}.
  \end{gather*}
  Transforming back to $I_j$ gives the result.
\end{proof}

For our spectral boundary layer mesh, the previous results give
\begin{gather*}
  \norm{\left( u - \mathcal{I}_p u \right)^{(k)}}_{0,\Omega_{BL}}
     \lesssim C_u \eps^{-1/2} h_j^{3/2-k} \E^{-\beta p}, \ \ k=0,1,2,
\end{gather*}
or equivalently, since $h_j = \kappa p \eps$,
\begin{gather}
  \label{res1p}
    (\kappa p \eps)^{-1} \norm{u-\mathcal{I}_p u}_{0,\Omega_{BL}}
     +  \abs{u-\mathcal{I}_p u}_{1,\Omega_{BL}}
     + \kappa p \eps \abs{u-\mathcal{I}_p u}_{2,\Omega_{BL}}
          \lesssim C_u \left(\kappa p  \right)^{1/2} \E^{-\beta p}. 
\end{gather}
Now, since \mbox{$\left(u- \mathcal{I}_p u \right)^{(k)}(0) = 0$},
\mbox{$k=0,1$}, we may write
\begin{gather*}
  \abs{\left(u- \mathcal{I}_p u\right)^{(k)}(x)}
    = \abs{\int_{0}^{x}\left(u- \mathcal{I}_p u\right)^{(k+1)}(t)\D t}, \ \ x \in \Omega_{BL}.
\end{gather*}
Then, by the Cauchy-Schwarz inequality and \eqref{res1p} we get for \mbox{$k=0,1$},
\begin{gather*}
  \abs{\left(u- \mathcal{I}_p u\right)^{(k)}(x)}
     \le (\kappa p \eps )^{1/2}
         \abs{\left(u- \mathcal{I}_p u\right)^{(k)}}_{1,\Omega_{BL}}
     \lesssim
       \begin{cases}
         C_u \eps^{-1/2} (\kappa p \eps)^{2} \E^{-\beta p}\;, & k=0, \\
         C_u \eps^{-1/2} \kappa p \eps   \E^{-\beta p} \;,    & k=1 .
       \end{cases}
\end{gather*}
Thus
\begin{gather*}
  (\kappa p \eps)^{-3/2} \norm{u-\mathcal{I}_p u}_{\infty,\Omega_{BL}}
     + (\kappa p \eps)^{-1/2} \norm{\left(u- \mathcal{I}_p u \right)'}_{\infty,\Omega_{BL}}
    \lesssim C_u (\kappa p )^{1/2} \E^{-\beta p}. 
\end{gather*}

\end{document}